\documentclass[preprint,10t]{elsarticle}

\usepackage{amssymb}
\usepackage{stmaryrd}
\usepackage{tipa}
\usepackage{amsfonts,amsmath,latexsym}
\usepackage{mathrsfs}
\usepackage{amsbsy}
\usepackage{indentfirst}
\usepackage{amsmath}
\usepackage{graphicx,psfrag,epsf}
\usepackage{enumerate}
\usepackage{natbib}
\usepackage{url} 
\usepackage{amsfonts}
\usepackage{mathrsfs}
\usepackage{algorithm}
\usepackage{algorithmic}
\usepackage{subfigure}
\usepackage{epsfig}
\usepackage{booktabs}
\usepackage{multirow}

\numberwithin{equation}{section}

\newtheorem{theorem}{\textbf{Theorem}}[section]
\newtheorem{lemma}{\textbf{Lemma}}

\newtheorem{cor}{\textbf{Corollary}}

\newtheorem{ass}{\textbf{Assumption}}

\newtheorem{remark}[theorem]{Remark}

\newcommand{\beqnar}{\begin{eqnarray*}}
\newcommand{\eeqnar}{\end{eqnarray*}}
\newcommand{\ba}{\begin{array}}
\newcommand{\ea}{\end{array}}

\newenvironment{proof}[1]{\begin{trivlist}\item {\it
\bf Proof.}\quad} {\qed\end{trivlist}}

\allowdisplaybreaks[4]

\journal{}

\begin{document}

\begin{frontmatter}



\title{ On Double Smoothed Volatility Estimation of Potentially Nonstationary Jump-Diffusion Model }
\author{Yuping Song\footnote{Corresponding author, email: songyuping@shnu.edu.cn }}

\address{ School of Finance and Business, Shanghai Normal
         University, Shanghai, 200234, P.R.China}

\begin{abstract}
In this paper, we present the double smoothed nonparametric approach
for infinitesimal conditional volatility of jump-diffusion model
based on high frequency data. Under certain minimal conditions, we
obtain the strong consistency and asymptotic normality for the
estimator as the time span $T \rightarrow \infty$ and the sample
interval $\Delta_{n} \rightarrow 0.$ The procedure and asymptotic
behavior can be applied for both null Harris recurrent and positive
Harris recurrent processes.
\end{abstract}


\end{frontmatter}

\section{Introduction}

Volatility is a very important variable in the research of financial
economics. Portfolio selection, original asset or derivative assets
option pricing and risk management depend on the accurate
measurement of volatility. In the past few decades, the estimation
for volatility has become one of the most active research fields in
empirical finance or time series econometrics.

Continuous-time models are widely used in economics and finance,
such as interest rate or an asset price, especially the
continuous-time diffusion processes, one of the most famous models
is the Black and Scholes asset pricing model. Bandi and Phillips
\cite{bp} proposed a double smoothed approach to the unknown
coefficients of potentially nonstationary or stationary scalar
diffusion models under high frequency data. Recently, an inordinate
amount of attention has been focused on the stochastic processes
with jumps, which can accommodate the impact of sudden and large
shocks to financial markets. Johannes \cite{joh} provided the
statistical and economic role of jumps in continuous-time interest
rate models. In this paper, we are focused on the following
time-homogeneous diffusion process with jumps $X = (X_{t})_{t \geq
0}$ such that
\begin{equation}
\label{jdm} dX_{t} = \left[\mu(X_{t-}) - \lambda(X_{t-})
\int_{Y}{c(X_{t-}, y)}\Pi(dy)\right]dt + \sigma(X_{t-})dW_{t} + d
J_{t},
\end{equation}
where $J_{t}$ is a compound Poisson process of the representation
\begin{equation}
d J_{t} = \int_{Y}{c(X_{t-}, y)}N(dt, dy),
\end{equation}
and $N(dt, dy)$ is a time-homogeneous Poisson counting measure with
independent increments. In addition, the functions $\mu(\cdot)$ and
$\sigma(\cdot)$ are the infinitesimal conditional drift and
variation due to the continuous diffusion process, $Y =
\mathbb{R}\setminus\{0\},W_{t}$ is a standard Brownian motion
independent of jump process $J_{t},$ $\lambda(\cdot)$ represents the
intensity measure and $c(\cdot, y)$ reflects the conditional
magnitude of a jump, where $y$ is a random variable with probability
measure $\Pi(dy).$ Furthermore, the coefficients $\mu(\cdot),$
$\sigma(\cdot),$ $\lambda(\cdot)$ and $c(\cdot, y)$ economically
represent the time trend referred as expected risk return and the
conditional variance of the return for an underlying asset,
respectively.

There are many statisticians and economists focused on state-domain
nonparametric estimation for volatility functions of jump-diffusion
models (\ref{jdm}). Bandi and Nguyen \cite{bn} and Johannes
\cite{joh} considered the kernel weighted version of instantaneous
volatility. Hanif \cite{hm1} \cite{hm2}, Hanif, Wang and Lin
\cite{hwl}, Lin and Wang \cite{lw2}, Wang, Zhang and Tang
\cite{wzt2}, Xu and Phillips \cite{xp} and so on improved the
properties of the estimators based on outstanding nonparametric
approaches. Moreover, Song \cite{sy}, Chen and Zhang \cite{cz}
considered the nonparametric volatility estimation of second-order
jump-diffusion model under integrated observations based on Wang and
Lin \cite{wl}. However, these work are only single-smoothing. In
this paper, we present the double smoothed nonparametric approach
for infinitesimal conditional volatility of jump-diffusion model
based on high frequency data. Double smoothed estimator can reduce
the asymptotic mean-squared error than the estimator proposed in
Bandi and Nguyen \cite{bn} for some chosen bandwidths, any level $x$
and any processes. Furthermore, for both null Harris recurrent and
positive Harris recurrent processes (\ref{jdm}), we obtain the
strong consistency and asymptotic normality for the estimator as the
time span $T \rightarrow \infty$ and the sample interval $\Delta_{n}
\rightarrow 0.$ Our result can availably solve the problem proposed
in the discussion part of Zhou \cite{zlk}.

The paper is organized as follows. The large sample properties of
the double smoothed volatility estimator are presented in Section 2.
Some technical lemmas and detailed proofs for the main theorem are
given in Section 3.

\section{ Technical assumptions and Large sample properties }

Model (\ref{jdm}) can be written in the integral form as
\begin{equation}
X_{t+\Delta} = X_{t} + \int_{t}^{t+\Delta} \mu(X_{t-})dt +
\int_{t}^{t+\Delta} \sigma(X_{t-})dW_{t} + \int_{t}^{t+\Delta}
\int_{Y}{c(X_{t-}, y)} \bar{\nu}(dt, dy),
\end{equation}
where $\bar{\nu}(dt, dy) := N(dt, dy) - \lambda(X_{t-}) \Pi(dy) dt$
is a compensated Poisson random measure. We can observe that
\begin{equation}
\int_{t}^{t+\Delta} \int_{Y}{c(X_{t-}, y)} \bar{\nu}(dt, dy) =
\int_{t}^{t+\Delta} dJ_{t} - \int_{t}^{t+\Delta}\lambda(X_{t-})
\int_{Y}{c(X_{t-}, y)}\Pi(dy)dt
\end{equation}
represents the conditional variation due to the discontinuous jumps
of the process $X_{t}$.

Due to the Markov properties of model (\ref{jdm}), we can build the
following infinitesimal conditional expectations as those in Bandi
and Nguyen \cite{bn}
\begin{align}
& \label{ice1} M^{2}(x) = \lim_{\Delta_{n} \rightarrow 0}
E\left[\frac{(X_{t + \Delta_{n}} - X_{t})^{2}}{\Delta_{n}}|X_{t} =
x\right]
 = \sigma^{2}(x) + \lambda(x)E_{Y}[c^{2}(x, y)], \\
& \label{ice2} M^{k}(x) = \lim_{\Delta_{n} \rightarrow 0}
E\left[\frac{(X_{t + \Delta_{n}} - X_{t})^{k}}{\Delta_{n}}|X_{t} =
x\right]
 = \lambda(x)E_{Y}[c^{k}(x, y)],
\end{align}
where $k > 2.$

Define
\begin{align*}
& t(i\Delta_{n})_{0} = \inf \{t \geq 0: |X_{t} - X_{i\Delta_{n}}| \leq \varepsilon_{n}\}, \\
& t(i\Delta_{n})_{j+1} = \inf \{t \geq t(i\Delta_{n})_{j} +
\Delta_{n}: |X_{t} - X_{i\Delta_{n}}| \leq \varepsilon_{n}\},\\
& m_{n}(i\Delta_{n}) = \sum_{j = 1}^{n} 1_{\{|X_{j\Delta_{n}} -
X_{i\Delta_{n}}| \leq \varepsilon_{n}\}}.
\end{align*}

For the given $\{X_{i\Delta_{n}} ; i= 1, 2, \cdots\}$, the double
smoothed estimator for volatility $M^{2}(x)$ based on the
infinitesimal conditional expectation (\ref{ice1}) is defined as
\begin{equation}
\label{dse} \hat{M}^{2}_{n}(x) = \frac{\sum_{i =
1}^{n}K_{h_{n}}(X_{i\Delta_{n}} - x) \frac{1}{m_{n}(i\Delta_{n})
\Delta_{n}} \sum_{j = 0}^{m_{n}(i\Delta_{n}) -
1}\left[X_{t(i\Delta_{n})_{j} + \Delta_{n}} -
X_{t(i\Delta_{n})_{j}}\right]^{2}}{\sum_{i =
1}^{n}K_{h_{n}}(X_{i\Delta_{n}} - x)},
\end{equation}
where $K_{h_{n}}( \cdot ) := \frac{1}{h_{n}} K(\frac{ \cdot
}{h_{n}})$ with the kernel function $K(\cdot)$ and $h_{n}$ is a
sequence of positive numbers, satisfies $h_n\to 0$ as $n\to \infty.$

The assumptions of this paper are listed below, which confirm the
large sample properties of the constructed estimators based on
(\ref{dse}). In what follows, denote $\mathscr{D} = (l, u)$ as the
admissible range of the process $X_{t}$ in model (\ref{jdm}).

\begin{ass}\label{a1}
i) ~ For each $n \in \mathbb{N},$ there exist a constant $C_{1}$ and
a function $\zeta_{n}$ : $Y \rightarrow \mathbb{R}_{+}$ with
$\int_{Y} \zeta_{n}^{2}(y) \Pi(dy) < \infty$ such that, for any $|x|
\leq n, |z| \leq n, y \in Y$,
$$|\mu(x) - \mu(z)| + |\sigma(x) - \sigma(z)| \leq L_{n}|x - z|,$$
$$ \lambda(x) \int_{Y} |c(x , y) - c(z , y)| \Pi(dy) \leq \zeta_{n}(y)|x - z|.$$

\noindent Moreover, for each $n \in \mathbb{N}$, there exist
$\zeta_{n}$ as above and $C_{2}$, such that for all $x \in
\mathbb{R}, y \in Y$,
$$|\mu(x)| + |\sigma(x)| \leq C_{2} (1 + |x|) , ~ \lambda(x) \int_{Y}|c(x , y)|\Pi(dy) \leq \zeta_{n}(y)(1 + |x|).$$

\noindent ii) ~ There exists a constant $C_{3}$ such that
$$\lambda(x) \int_{Y} |c(x, y)|^{\alpha} \Pi(dy) \leq C_{3} (1 + |x|^{\alpha}).$$
for a fixed $\alpha > 2$ and $\forall x \in \mathscr{D}.$
\medskip

\noindent iii) ~ The functions $\sigma^{2}( \cdot ),$ $\lambda(
\cdot )$ and $c( \cdot, y )$ are at least twice continuously
differentiable. $\lambda( x ) \geq 0$ and $\sigma^{2}( x ) \geq 0$
for $\forall x \in \mathscr{D}.$
\end{ass}

\begin{remark}
\label{ar1} \rm{This assumption guarantees the existence and
uniqueness of a c\`{a}dl\`{a}g strong solution to stochastic
differential equation $X_{t}$ in (\ref{jdm}), see Jacod and Shiryaev
\cite{js}. For instance, Bandi and Nguyen \cite{bn}, Shimizu and
Yoshida \cite{syo} imposed similar conditions on the coefficients of
the underlying stochastic differential equation. }
\end{remark}

\begin{ass} \label{a2}
The process $X=\{X_{t}\}_{\ge 0}$ in model (\ref{jdm}) is Harris
recurrent.
\end{ass}

\begin{ass} \label{a3}
The process $X=\{X_{t}\}_{\ge 0}$ in model (\ref{jdm}) is positive
Harris recurrent.
\end{ass}

\begin{remark}
\label{ar2} \rm{The Assumption \ref{a2} guarantees the existence of
a unique invariant measure $s(x),$ that is, $s(A) =
\int_{\mathscr{D}}P(X_{t}^{(x)} \in A)s(dx) ~~ \forall A \in
\mathfrak{B}(\mathscr{D}).$ The Assumption \ref{a3} implies that the
process $X_{t}$ has a time-invariant probability measure given by
$p(dx) = \frac{s(x)}{s(\mathscr{D})}.$ The positive Harris recurrent
condition means that the process becomes stationary at any initial
level $x \in \mathscr{D}.$ Furthermore, as discussed in Bandi and
Phillips \cite{bp}, the stationary probability measure can increase
the asymptotic rate of convergence for underlying estimators.}
\end{remark}

\begin{ass} \label{a4}
The kernel $K$($\cdot$) : $\mathbb{R} \rightarrow \mathbb{R}^{+}$ is
a continuously differentiable, bounded and symmetric function
satisfying:
$$\int K(u) du = 1,~\int K^{'}(u) du < \infty,~K_{i}^{j} := \int K^{i}(u) u^{j} du < \infty.$$
\end{ass}

\begin{remark}
\label{ar3} \rm{In fact, any density function can be considered as a
kernel, moreover even unnecessary positive functions can be used.
For simplification, we only consider positive and symmetrical
kernels used widely. It is well known both empirically and
theoretically that the choice of kernel functions is not very
important to the kernel estimator, see Gasser and M\"{u}ller
\cite{gm}.}
\end{remark}

\begin{ass} \label{a5}
$T \rightarrow \infty,~\Delta_{n} \rightarrow 0,~h_{n}\rightarrow
0,~\frac{\bar{L}_{X}(T,
x)}{h_{n}}(\Delta_{n}\log(\frac{1}{\Delta_{n}}))^{\frac{1}{2}}
\rightarrow 0,~~\varepsilon_{n}\bar{L}_{X}(T, x) \rightarrow
\infty,$ $\frac{\bar{L}_{X}(T,
x)}{\varepsilon_{n}}(\Delta_{n}\log(\frac{1}{\Delta_{n}}))^{\frac{1}{2}}
\rightarrow 0 ~as~ n \rightarrow \infty.$
\end{ass}

\begin{remark}
\label{ar4}\rm{ The relationship between $h_{n}$ and $\Delta_{n}$ is
similar as that in Bandi and Nguyen \cite{bn}.}
\end{remark}

We have the following asymptotic results for the double smoothed
estimators such as (\ref{dse}) based on the assumptions above.

\begin{theorem}
\label{mr} Under Assumptions \ref{a1}, \ref{a2}, \ref{a4}, \ref{a5},
as $n \rightarrow \infty,$ we have

(i)~~$\hat{M}^{2}_{n}(x) \stackrel{a.s.} \longrightarrow M^{2}(x).$

(ii)~Furthermore, if $h_{n} = o(\varepsilon_{n})$ and
$\varepsilon^{n}_{5} \bar{L}_{X}(T, x) = O_{a.s.}(1),$ then
\begin{equation*}
\sqrt{\varepsilon_{n} \hat{\bar{L}}_{X}(T, x)} (\hat{M}^{2}_{n}(x) -
M^{2}(x) - \varepsilon_{n}^{2} \ast \Gamma_{M^{2}}) \Rightarrow
N\left(0, \frac{1}{2} M^{4}(x)\right),
\end{equation*}
where $\Gamma_{M^{2}} = \frac{1}{3} \left[\frac{1}{2}(M^{2})^{''}(x)
+ (M^{2})^{'}(x)\frac{s^{'}(x)}{s(x)} \right],$

(iii)~If $h_{n} = O(\varepsilon_{n})$ with
$\frac{h_{n}}{\varepsilon_{n}} = \phi$ and $\varepsilon^{n}_{5}
\bar{L}_{X}(T, x) = O_{a.s.}(1),$ then
\begin{equation*}
\sqrt{\varepsilon_{n} \hat{\bar{L}}_{X}(T, x)} (\hat{M}^{2}_{n}(x) -
M^{2}(x) - \varepsilon_{n}^{2} \ast \Gamma_{M^{2}}^{\phi})
\Rightarrow N\left(0, \frac{1}{2} \theta_{\phi} M^{4}(x)\right),
\end{equation*}
where $\Gamma_{M^{2}}^{\phi} = \left(K_{2}\phi^{2} +
\frac{1}{3}\right)\left[\frac{1}{2}(M^{2})^{''}(x) +
(M^{2})^{'}(x)\frac{s^{'}(x)}{s(x)}\right]$ and

\noindent $\theta_{\phi} = \frac{1}{2}\int_{- \infty}^{+ \infty} dz
\int_{(z - 1)/\phi}^{(z + 1)/\phi} da \int_{(z - 1)/\phi}^{(z +
1)/\phi} de K\left(a\right) K\left(e\right),$
\end{theorem}
where $\hat{\bar{L}}_{X}(T, x)$ is defined as that in lemma
\ref{lk}.
\bigskip

Under Assumption \ref{a3}, the local time $\bar{L}^{\oplus}_{X}(t,
a)$ increases consistently with $T$ up to multiplication by a
constant as
\begin{equation}
\label{ltt} \frac{\bar{L}^{\oplus}_{X}(t, a)}{T}
\stackrel{a.s.}{\longrightarrow} p(x),~\forall x \in \mathscr{D},
\end{equation},
where $\bar{L}^{\oplus}_{X}(t, a)$ is mentioned in equation
(\ref{lt6}) below. Based on the equation (\ref{ltt}), we can obtain
the following corollary.

\begin{cor}
\label{mr2} [\textbf{Stationary~Case}] Under Assumptions \ref{a1},
\ref{a3}, \ref{a4}, \ref{a5}, as $n \rightarrow \infty,$ we have

(i)~~$\hat{M}^{2}_{n}(x) \stackrel{P} \longrightarrow M^{2}(x).$

(ii)~Furthermore, if $h_{n} = o(\varepsilon_{n})$ and $n \Delta_{n}
\varepsilon^{n}_{5} = O_{a.s.}(1),$ then
\begin{equation*}
\sqrt{n \Delta_{n} \varepsilon_{n}} (\hat{M}^{2}_{n}(x) - M^{2}(x) -
\varepsilon_{n}^{2} \ast \Gamma_{M^{2}}) \Rightarrow N\left(0,
\frac{1}{2} \frac{M^{4}(x)}{p(x)} \right),
\end{equation*}
where $\Gamma_{M^{2}} = \frac{1}{3} \left[\frac{1}{2}(M^{2})^{''}(x)
+ (M^{2})^{'}(x)\frac{s^{'}(x)}{s(x)} \right],$

(iii)~If $h_{n} = O(\varepsilon_{n})$ with
$\frac{h_{n}}{\varepsilon_{n}} = \phi$ and $n \Delta_{n}
\varepsilon^{n}_{5} = O_{a.s.}(1),$ then
\begin{equation*}
\sqrt{n \Delta_{n} \varepsilon_{n}} (\hat{M}^{2}_{n}(x) - M^{2}(x) -
\varepsilon_{n}^{2} \ast \Gamma_{M^{2}}^{\phi}) \Rightarrow
N\left(0, \frac{1}{2} \theta_{\phi} \frac{M^{4}(x)}{p(x)}\right),
\end{equation*}
where $\Gamma_{M^{2}}^{\phi} = \left(K_{2}\phi^{2} +
\frac{1}{3}\right)\left[\frac{1}{2}(M^{2})^{''}(x) +
(M^{2})^{'}(x)\frac{s^{'}(x)}{s(x)}\right]$ and

\noindent $\theta_{\phi} = \frac{1}{2}\int_{- \infty}^{+ \infty} dz
\int_{(z - 1)/\phi}^{(z + 1)/\phi} da \int_{(z - 1)/\phi}^{(z +
1)/\phi} de K\left(a\right) K\left(e\right).$
\end{cor}

\begin{remark}
\label{mrr1} \rm{ In contrary to the scalar diffusion model without
jumps in Bandi and Phillips \cite{bp}, the rate of convergence of
the second infinitesimal moment estimator is same as the first
infinitesimal moment estimator. Apparently, this is due to the
presence of discontinuous breaks that have an equal impact on all
the functional estimates. As Johannes \cite{joh} pointed out, for
the conditional variance of interest rate changes, not only
diffusion play a certain role, but also jumps account for more than
half at lower interest level rates, almost two-thirds at higher
interest level rates, which dominate the conditional volatility of
interest rate changes. Thus, it is extremely important to estimate
the conditional variance as $M^{2}(x)$ not only the diffusion part
$\sigma^{2}(x)$, which reflects the fluctuation of the return of the
underlying asset. Nonparametric estimation to identify the diffusion
coefficient $\sigma^{2}(x)$, the jump intensity $\lambda(x)$ and the
jump sizes $c(x, y)$ for model (\ref{jdm}) is not our objective in
this paper and thus it is less of a concern here and left for the
future research.}
\end{remark}

\begin{remark}
\label{mrr2} \rm{There are many statisticians and economists focused
on state-domain nonparametric estimation for volatility functions of
diffusion models with jumps. Bandi and Nguyen \cite{bn} and Johannes
\cite{joh} considered the kernel weighted version of instantaneous
volatility combinated with the combination of power variation when
the price process follows scalar diffusion model with jumps as
(\ref{jdm}). They established the following asymptotic normality for
the estimator $\hat{M}^{2}_{bn}(x)$ of unknown quantity $M^{2}(x)$
with the Assumptions \ref{a1}, \ref{a2}, \ref{a4} and \ref{a5}, that
is,
\begin{equation*}
\sqrt{h_{n} \hat{\bar{L}}_{X}(T, x)} (\hat{M}^{2}_{bn}(x) - M^{2}(x)
- \Gamma_{M^{2}}^{bn}) \Rightarrow N\left(0, K_{2}M^{4}(x)\right),
\end{equation*}
where $\hat{M}^{2}_{bn}(x) = \frac{\sum_{i =
1}^{n}K_{h_{n}}(X_{i\Delta_{n}} - x) \left(X_{(i + 1) \Delta_{n}} -
X_{i\Delta_{n}}\right)^{2} }{\Delta_{n} \sum_{i =
1}^{n}K_{h_{n}}(X_{i\Delta_{n}} - x)}$ and

\noindent $\Gamma_{M^{2}}^{bn} = h_{n}^{2}
\left[\frac{1}{2}(M^{2})^{''}(x) +
(M^{2})^{'}(x)\frac{s^{'}(x)}{s(x)} \right].$

As mentioned in Bandi and Phillips \cite{bp}, the estimator
$\hat{M}^{2}_{bn}(x)$ proposed in Bandi and Nguyen \cite{bn} is the
same as the double smoothed estimator $\hat{M}^{2}(x)$ conducted as
(\ref{dse}) asymptotically if $h_{n} = o(\varepsilon_{n}).$
Moreover, as discussed in Bandi and Phillips \cite{bp}, if
$\frac{h_{n}}{\varepsilon_{n}} = \phi,$ double smoothed estimator
can reduce the asymptotic mean-squared error than the estimator
$\hat{M}^{2}_{bn}(x)$ above for some chosen bandwidth $h_{n}$, any
level $x$ and any processes. }
\end{remark}

\begin{remark}
\label{mrr3} \rm{It is very important to consider the choice of the
bandwidth in nonparametric estimation. Here we will select the
optimal bandwidth $h_{n}$ based on the mean squared error (MSE) and
the asymptotic theory in Theorem \ref{mr}. The optimal smoothing
parameter $h_{n}$ for double smoothed estimator of $M^{2}(x)$ is
given that
\begin{equation*}
h_{n,opt} = \phi \ast \left(\frac{1}{\bar{L}_{X}(T, x)} \cdot
\frac{\frac{1}{2} \theta_{\phi} M^{4}(x) }{\left(K_{2}\phi^{2} +
\frac{1}{3}\right)^{2}\left[\frac{1}{2}(M^{2})^{''}(x) +
(M^{2})^{'}(x)\frac{s^{'}(x)}{s(x)}\right]^{2} }
\right)^{\frac{1}{5}}
 = O_{p}\left(\frac{1}{\bar{L}_{X}(T, x)}\right)^{\frac{1}{5}},\end{equation*}
which differs from the continuous case in Bandi and Phillips
\cite{bp} with
$h_{n,opt}=O_{p}\left(\frac{\Delta_{n}}{\bar{L}_{X}(T,
x)}\right)^{\frac{1}{5}}.$ Furthermore, one can discuss the optimal
bandwidth for double smoothed volatility estimator of jump-diffusion
model based on Wang and Zhou \cite{wz}, which will be under
consideration in the future study. If the smoothing parameter $h_{n}
= O((\bar{L}_{X}(T, x))^{-1/5}),$ the normal confidence interval for
$M^{2}(x)$ using double smoothed estimators at the significance
level $100(1-\alpha)\%$ are constructed as follows,
\begin{align*}
I_{\mu,\alpha} = & \Bigg[\hat{M}^{2}_{n}(x) - \varepsilon_{n}^{2}
\cdot \hat{\Gamma}_{\hat{M}^{2}}^{\phi} - z_{1-\alpha/2} \cdot
\frac{1}{\sqrt{\varepsilon_{n} \hat{\bar{L}}_{X}(T, x)}} \cdot
\sqrt{\frac{1}{2}
\theta_{\phi} \hat{M}^{4}(x)},\\
& \hat{M}^{2}_{n}(x) - \varepsilon_{n}^{2} \cdot
\hat{\Gamma}_{\hat{M}^{2}}^{\phi} + z_{1-\alpha/2} \cdot
\frac{1}{\sqrt{\varepsilon_{n} \hat{\bar{L}}_{X}(T, x)}} \cdot
\sqrt{\frac{1}{2} \theta_{\phi} \hat{M}^{4}(x)} \Bigg],
\end{align*}
where $z_{1-\alpha/2}$ is the inverse CDF for the standard normal
distribution evaluated at $1 - \alpha/2.$ To facilitate statistical
inference for $M_{2}(x)$ based on Theorem \ref{mr}, we need to
conduct consistent estimators for the unknown quantities $M^{4}(x)$
in the normal approximation. Based on the infinitesimal moments
condition (\ref{ice2}), the double smoothed estimator for $M^{4}(x)$
is conducted as
\begin{equation}
\label{dse} \hat{M}^{4}_{n}(x) = \frac{\sum_{i =
1}^{n}K_{h_{n}}(X_{i\Delta_{n}} - x) \frac{1}{m_{n}(i\Delta_{n})
\Delta_{n}} \sum_{j = 0}^{m_{n}(i\Delta_{n}) -
1}\left[X_{t(i\Delta_{n})_{j} + \Delta_{n}} -
X_{t(i\Delta_{n})_{j}}\right]^{4}}{\sum_{i =
1}^{n}K_{h_{n}}(X_{i\Delta_{n}} - x)}.
\end{equation}
The consistency and asymptotic normality for $\hat{M}^{4}_{n}(x)$
can be done with the similar approach as $\hat{M}^{2}_{n}(x),$ which
goes beyond the scope here and will be under consideration in the
future study. }
\end{remark}

\section{ Detailed Proof }
In this section, we first present some technical lemmas and the
proofs for the main theorems.

\subsection{ Some Technical Lemmas with Proofs}

\begin{lemma}
\label{lot} (Bandi and Nguyen \cite{bn}) Let $X$ be a semimartingale
with local time $L_{X}(\cdot, a)_{a \in \mathscr{D}}$ and $f$ be a
bounded Borel measurable function, we have
\begin{equation} \label{ot}
\int_{0}^{t} g(X_{s-}) d [X]_{s}^{c} = \int_{- \infty}^{\infty}
L_{X}(t, a) g(a) da,~~~~~~~~~~a.s.
\end{equation}
where $[X]_{s}^{c}$ denotes the continuous part of the quadratic
variation of $X.$
\end{lemma}

\begin{lemma}
\label{llt} (Bandi and Nguyen \cite{bn}) Let $X$ be a semimartingale
satisfying $\sum_{0 < s \leq t} |\Delta X_{s}| < \infty$ a.s.
$\forall t.$ Then, $\forall (t, a)$ we have
\begin{equation} \label{lt1}
L_{X}(t, a+) = L_{X}(t, a) = \lim_{\varepsilon \rightarrow 0}
\frac{1}{\varepsilon} \int_{0}^{t} 1_{(a \leq X_{s} \leq a +
\varepsilon)}d[X]_{s}^{c},~~~~~~~~~~a.s.
\end{equation}
and
\begin{equation} \label{lt2}
L_{X}(t, a-) = \lim_{\varepsilon \rightarrow 0}
\frac{1}{\varepsilon} \int_{0}^{t} 1_{(a - \varepsilon \leq X_{s}
\leq a) }d[X]_{s}^{c},~~~~~~~~~~a.s.
\end{equation}
Also,
\begin{equation} \label{lt3}
L^{\oplus}_{X}(t, a) = \frac{L_{X}(t, a) + L_{X}(t, a-)}{2} =
\lim_{\varepsilon \rightarrow 0} \frac{1}{2\varepsilon} \int_{0}^{t}
1_{(| X_{s} - a| \leq \varepsilon) }d[X]_{s}^{c},~~~~~~~~~~a.s.
\end{equation}
\end{lemma}

\begin{remark}
\label{rlt}\rm{ We may employ the following versions of local time
in what follows.
\begin{equation} \label{lt4}
\bar{L}_{X}(t, a+) = \bar{L}_{X}(t, a) =
\frac{1}{\sigma^{2}(a)}\lim_{\varepsilon \rightarrow 0}
\frac{1}{\varepsilon} \int_{0}^{t} 1_{(a \leq X_{s} \leq a +
\varepsilon)} \sigma^{2}(X_{s}) ds,~~~~~~~~~~a.s.
\end{equation}
and
\begin{equation} \label{lt5}
\bar{L}_{X}(t, a-) = \frac{1}{\sigma^{2}(a)}\lim_{\varepsilon
\rightarrow 0} \frac{1}{\varepsilon} \int_{0}^{t} 1_{(a -
\varepsilon \leq X_{s} \leq a) } \sigma^{2}(X_{s}) ds,~~~~~~~~~~a.s.
\end{equation}
Also,
\begin{equation} \label{lt6}
\bar{L}^{\oplus}_{X}(t, a) = \frac{\bar{L}_{X}(t, a) +
\bar{L}_{X}(t, a-)}{2} = \frac{1}{\sigma^{2}(a)}\lim_{\varepsilon
\rightarrow 0} \frac{1}{2\varepsilon} \int_{0}^{t} 1_{(| X_{s} - a|
\leq \varepsilon) } \sigma^{2}(X_{s}) ds,~~~~~~~~~~a.s.
\end{equation}
}
\end{remark}

\begin{lemma}
\label{lk} (Bandi and Nguyen \cite{bn}) Under Assumptions \ref{a1} -
\ref{a5}, we have
\begin{equation} \label{kk}
\hat{\bar{L}}_{X}(T, x) =
\frac{\Delta_{n}}{h_{n}}\sum_{i=1}^{n}K\left(\frac{X_{i\Delta_{n}} -
x}{h_{n}}\right) \stackrel{\mathrm{a.s.}}{\longrightarrow}
\bar{L}^{\oplus}_{X}(T, x).
\end{equation}
\end{lemma}

\subsection{ The proof of Theorem \ref{mr} }
\begin{proof}
\medskip

\noindent {\textbf{Strong~Consistency:}}
\begin{eqnarray*}
\hat{M}^{2}_{n}(x) & = & \frac{\sum_{i =
1}^{n}K_{h_{n}}(X_{i\Delta_{n}} - x) \frac{1}{m_{n}(i\Delta_{n})
\Delta_{n}} \sum_{j = 0}^{m_{n}(i\Delta_{n}) -
1}\left(\left[X_{t(i\Delta_{n})_{j} + \Delta_{n}} -
X_{t(i\Delta_{n})_{j}}\right]^{2} -
M^{2}(X_{i\Delta_{n}})\right)}{\sum_{i =
1}^{n}K_{h_{n}}(X_{i\Delta_{n}} - x)}\\
& ~ & - \frac{\sum_{i = 1}^{n}K_{h_{n}}(X_{i\Delta_{n}} - x)
M^{2}(X_{i\Delta_{n}}) }{\sum_{i = 1}^{n}K_{h_{n}}(X_{i\Delta_{n}} -
x)}\\
& := & A_{1n} - A_{2n}.
\end{eqnarray*}
Now, define
$$\delta_{n, T}  = \max_{i \leq n} \sup_{i\Delta_{n} \leq s \leq (i+1)\Delta_{n}}\left|X_{s-} - X_{i\Delta_{n}}\right|.$$
As is shown in Bandi and Nguyen \cite{bn} that
\begin{equation}
\label{holder} \varlimsup_{n \rightarrow \infty} {\frac{\delta_{n,
T}}{\left(\Delta_{n} \log(1/\Delta_{n})\right)^{1/2}}} = C_{1}
~~~~~~ a.s.
\end{equation}
for some constant $C_{1},$ which implies that $\delta_{n, T} =
o_{a.s.}(1).$

As for $A_{1n},$ based on equation (\ref{holder}) and the quotient
limit theorem for Harris recurrent Markov processes, we can obtain
that
\begin{eqnarray*}
A_{1n} & = & \frac{\int_{0}^{T}K_{h_{n}}(X_{s-} - x) M^{2}(X_{s-})ds
+ O_{a.s.}\left(\frac{\bar{L}_{X}(T, x)}{h_{n}}\left(\Delta_{n} \log
(\frac{1}{\Delta_{n}})\right)^{1/2}\right)
}{\int_{0}^{T}K_{h_{n}}(X_{s-} - x) ds +
O_{a.s.}\left(\frac{\bar{L}_{X}(T, x)}{h_{n}}\left(\Delta_{n} \log
(\frac{1}{\Delta_{n}})\right)^{1/2}\right)}\\
& = & \frac{M^{2}(x) s(x) + o_{a.s.}(1)}{s(x) + o_{a.s.}(1)} +
o_{a.s.}(1)\\
& = & M^{2}(x) = \sigma^{2}(x) + \lambda(x) E_{Y}[c^{2}(x, y)].
\end{eqnarray*}
For the term $A_{2n},$ it is sufficient to prove that
\begin{equation}
\label{a2n} \frac{1}{m_{n}(i\Delta_{n}) \Delta_{n}} \sum_{j =
0}^{m_{n}(i\Delta_{n}) - 1}\left[X_{t(i\Delta_{n})_{j} + \Delta_{n}}
- X_{t(i\Delta_{n})_{j}}\right]^{2} - M^{2}(X_{i\Delta_{n}})
\stackrel{a.s.}{\longrightarrow} 0.
\end{equation}
Using It\^{o} formula to the jump-diffusion setting shown in Protter
\cite{pr}, we can write
\begin{eqnarray*}
& ~ & (X_{t(i\Delta_{n})_{j} + \Delta_{n}} - X_{t(i\Delta_{n})_{j}})^{2}\\
 & = & 2\int_{t(i\Delta_{n})_{j}}^{t(i\Delta_{n})_{j} + \Delta_{n}}(X_{s-} -
 X_{t(i\Delta_{n})_{j}})\mu(X_{s-})ds + 2\int_{t(i\Delta_{n})_{j}}^{t(i\Delta_{n})_{j} + \Delta_{n}}(X_{s-} -
 X_{t(i\Delta_{n})_{j}})\sigma(X_{s-})dW_{s}\\
 & ~ & + 2\int_{t(i\Delta_{n})_{j}}^{t(i\Delta_{n})_{j} + \Delta_{n}}(X_{s-} -
 X_{t(i\Delta_{n})_{j}}) \int_{Y}c(X_{s-}, y) \bar{\nu}(ds, dy) +
 \int_{t(i\Delta_{n})_{j}}^{t(i\Delta_{n})_{j} + \Delta_{n}} M^{2}(X_{s-})ds\\
 & ~ & +
 \int_{t(i\Delta_{n})_{j}}^{t(i\Delta_{n})_{j} + \Delta_{n}} \int_{Y}c^{2}(X_{s-}, y) \bar{\nu}(ds,
 dy),
\end{eqnarray*}
which implies that (\ref{a2n}) can be divided into five parts as
\begin{eqnarray*}
& ~ & \frac{1}{m_{n}(i\Delta_{n}) \Delta_{n}} \sum_{j =
0}^{m_{n}(i\Delta_{n}) - 1}\left[X_{t(i\Delta_{n})_{j} + \Delta_{n}}
- X_{t(i\Delta_{n})_{j}}\right]^{2} - M^{2}(X_{i\Delta_{n}})\\
& = & \frac{1}{m_{n}(i\Delta_{n}) \Delta_{n}} \sum_{j =
0}^{m_{n}(i\Delta_{n}) - 1}
\int_{t(i\Delta_{n})_{j}}^{t(i\Delta_{n})_{j} + \Delta_{n}} \left(
M^{2}(X_{s-}) -
M^{2}(X_{i\Delta_{n}})\right)ds + \\
& + & \frac{2}{m_{n}(i\Delta_{n}) \Delta_{n}} \sum_{j =
0}^{m_{n}(i\Delta_{n}) -
1}\int_{t(i\Delta_{n})_{j}}^{t(i\Delta_{n})_{j} + \Delta_{n}}(X_{s-}
- X_{t(i\Delta_{n})_{j}})\mu(X_{s-})ds +\\
& + & \frac{2}{m_{n}(i\Delta_{n}) \Delta_{n}} \sum_{j =
0}^{m_{n}(i\Delta_{n}) -
1}\int_{t(i\Delta_{n})_{j}}^{t(i\Delta_{n})_{j} + \Delta_{n}}(X_{s-}
- X_{t(i\Delta_{n})_{j}})\sigma(X_{s-})dW_{s}\\
& + & \frac{2}{m_{n}(i\Delta_{n}) \Delta_{n}} \sum_{j =
0}^{m_{n}(i\Delta_{n}) - 1}
\int_{t(i\Delta_{n})_{j}}^{t(i\Delta_{n})_{j} + \Delta_{n}}(X_{s-} -
 X_{t(i\Delta_{n})_{j}}) \int_{Y}c(X_{s-}, y) \bar{\nu}(ds, dy)\\
& + & \frac{1}{m_{n}(i\Delta_{n}) \Delta_{n}} \sum_{j =
0}^{m_{n}(i\Delta_{n}) -
1} \int_{t(i\Delta_{n})_{j}}^{t(i\Delta_{n})_{j} + \Delta_{n}} \int_{Y}c^{2}(X_{s-}, y) \bar{\nu}(ds, dy)\\
& := & A_{21n} + A_{22n} + A_{23n} + A_{24n} + A_{25n}.
\end{eqnarray*}
For instance, we can write $A_{25n}$ as
\begin{eqnarray*}
& ~ & \frac{1}{m_{n}(i\Delta_{n}) \Delta_{n}} \sum_{j =
0}^{m_{n}(i\Delta_{n}) - 1}
\int_{t(i\Delta_{n})_{j}}^{t(i\Delta_{n})_{j} + \Delta_{n}}
\int_{Y}c^{2}(X_{s-}, y) \bar{\nu}(ds, dy)\\
& = & \frac{\sum_{j = 1}^{n}1_{\{|X_{j\Delta_{n}} - X_{i\Delta_{n}}|
\leq \varepsilon_{n}\}}\int_{j\Delta_{n}}^{(j + 1)\Delta_{n}}
\int_{Y} c^{2}(X_{s-}, y) \bar{\nu}(ds, dy) }{\Delta_{n} \sum_{j =
1}^{n}1_{\{|X_{j\Delta_{n}} - X_{i\Delta_{n}}| \leq
\varepsilon_{n}\}}}.
\end{eqnarray*}
$A_{23n}, A_{24n}, A_{25n}$ are sample averages of martingale
difference sequences, which converge to zero a.s.

Due to the locally boundedness of $\mu(\cdot),$ the term $A_{22n}$
can be done similarly as $A_{21n},$ here we only prove $A_{21n}
\stackrel{a.s.}{\longrightarrow} 0$ for simplicity.

By the mean-value theorem, the equation (\ref{holder}) and the
locally boundedness of $(M^{2})^{'}(\cdot),$ we can obtain
\begin{eqnarray*}
A_{21n} & = & \frac{1}{m_{n}(i\Delta_{n}) \Delta_{n}} \sum_{j =
0}^{m_{n}(i\Delta_{n}) - 1}
\int_{t(i\Delta_{n})_{j}}^{t(i\Delta_{n})_{j} + \Delta_{n}} \left(
M^{2}(X_{s-}) - M^{2}(X_{i\Delta_{n}})\right) ds\\
& = & \frac{1}{m_{n}(i\Delta_{n}) \Delta_{n}} \sum_{j =
0}^{m_{n}(i\Delta_{n}) - 1}
\int_{t(i\Delta_{n})_{j}}^{t(i\Delta_{n})_{j} + \Delta_{n}} \left(
M^{2}(X_{s-}) - M^{2}(X_{j\Delta_{n}})\right) ds\\
& + & \frac{1}{m_{n}(i\Delta_{n}) \Delta_{n}} \sum_{j =
0}^{m_{n}(i\Delta_{n}) - 1}
\int_{t(i\Delta_{n})_{j}}^{t(i\Delta_{n})_{j} + \Delta_{n}} \left(
M^{2}(X_{j\Delta_{n}}) - M^{2}(X_{i\Delta_{n}})\right) ds\\
& = & \frac{\sum_{j = 1}^{n}1_{\{|X_{j\Delta_{n}} - X_{i\Delta_{n}}|
\leq \varepsilon_{n}\}}\int_{j\Delta_{n}}^{(j + 1)\Delta_{n}} \left(
M^{2}(X_{s-}) - M^{2}(X_{j\Delta_{n}})\right)ds }{\Delta_{n} \sum_{j
= 1}^{n}1_{\{|X_{j\Delta_{n}} - X_{i\Delta_{n}}| \leq
\varepsilon_{n}\}}}\\
& + & \frac{\sum_{j = 1}^{n}1_{\{|X_{j\Delta_{n}} - X_{i\Delta_{n}}|
\leq \varepsilon_{n}\}}\int_{j\Delta_{n}}^{(j + 1)\Delta_{n}} \left(
M^{2}(X_{j\Delta_{n}}) - M^{2}(X_{i\Delta_{n}})\right)ds
}{\Delta_{n} \sum_{j = 1}^{n}1_{\{|X_{j\Delta_{n}} -
X_{i\Delta_{n}}| \leq
\varepsilon_{n}\}}}\\
& = & \frac{\sum_{j = 1}^{n}1_{\{|X_{j\Delta_{n}} - X_{i\Delta_{n}}|
\leq \varepsilon_{n}\}}\int_{j\Delta_{n}}^{(j + 1)\Delta_{n}}
(M^{2})^{'}(\xi_{n,1}) \left( X_{s-} - X_{j\Delta_{n}} \right)ds
}{\Delta_{n} \sum_{j = 1}^{n}1_{\{|X_{j\Delta_{n}} -
X_{i\Delta_{n}}| \leq
\varepsilon_{n}\}}}\\
& + & \frac{\sum_{j = 1}^{n}1_{\{|X_{j\Delta_{n}} - X_{i\Delta_{n}}|
\leq \varepsilon_{n}\}}\int_{j\Delta_{n}}^{(j + 1)\Delta_{n}}
(M^{2})^{'}(\xi_{n,2}) \left( X_{j\Delta_{n}} - X_{i\Delta_{n}}
\right)ds }{\Delta_{n} \sum_{j = 1}^{n}1_{\{|X_{j\Delta_{n}} -
X_{i\Delta_{n}}| \leq
\varepsilon_{n}\}}}\\
& \simeq & O_{a.s.}\left(\Delta_{n} \log (1/\Delta_{n})\right)^{1/2}
+ C \varepsilon_{n} \stackrel{a.s.}{\longrightarrow} 0,
\end{eqnarray*}
where $\xi_{n,1}$ lies between $X_{s-}$ and $X_{j\Delta_{n}},$
$\xi_{n,2}$ lies between $X_{1\Delta_{n}}$ and $X_{j\Delta_{n}}.$

We have proved that
\begin{equation}
\hat{M}^{2}_{n}(x) \stackrel{a.s.}{\longrightarrow} M^{2}(x) =
\sigma^{2}(x) + \lambda(x) E_{Y}[c^{2}(x, y)].
\end{equation}
\medskip
\noindent {\textbf{Asymptotic~Normality:}}
\begin{eqnarray*}
& ~ & \hat{M}^{2}_{n}(x) - M^{2}(x)\\
& = & \frac{ \Delta_{n} \sum_{i = 1}^{n}K_{h_{n}}(X_{i\Delta_{n}} -
x) \frac{1}{m_{n}(i\Delta_{n}) \Delta_{n}} \sum_{j =
0}^{m_{n}(i\Delta_{n}) - 1}\left[X_{t(i\Delta_{n})_{j} + \Delta_{n}}
- X_{t(i\Delta_{n})_{j}}\right]^{2}}{ \Delta_{n} \sum_{i
= 1}^{n}K_{h_{n}}(X_{i\Delta_{n}} - x)} - M^{2}(x)\\
& = & \frac{ \Delta_{n} \sum_{i = 1}^{n}K_{h_{n}}(X_{i\Delta_{n}} -
x) \frac{1}{m_{n}(i\Delta_{n}) \Delta_{n}} \sum_{j =
0}^{m_{n}(i\Delta_{n}) - 1}\left[X_{t(i\Delta_{n})_{j} + \Delta_{n}}
- X_{t(i\Delta_{n})_{j}}\right]^{2}}{ \Delta_{n} \sum_{i =
1}^{n}K_{h_{n}}(X_{i\Delta_{n}} - x)}\\
& ~ & - \frac{ \Delta_{n} \sum_{i = 1}^{n}K_{h_{n}}(X_{i\Delta_{n}}
- x) M^{2}(X_{i\Delta_{n}}) }{ \Delta_{n} \sum_{i =
1}^{n}K_{h_{n}}(X_{i\Delta_{n}} - x)}\\
& ~ & + \frac{ \Delta_{n} \sum_{i = 1}^{n}K_{h_{n}}(X_{i\Delta_{n}}
- x) M^{2}(X_{i\Delta_{n}}) }{ \Delta_{n} \sum_{i =
1}^{n}K_{h_{n}}(X_{i\Delta_{n}} - x)} - \frac{ \Delta_{n} \sum_{i =
1}^{n}K_{h_{n}}(X_{i\Delta_{n}} - x) M^{2}(x) }{ \Delta_{n} \sum_{i
= 1}^{n}K_{h_{n}}(X_{i\Delta_{n}} - x)}\\
& = & V + B.
\end{eqnarray*}
As for the bias term $B,$ based on equation (\ref{holder}) and the
quotient limit theorem for Harris recurrent Markov processes, we can
obtain that
\begin{eqnarray*}
B & = & \frac{ \frac{\Delta_{n}}{h_{n}} \sum_{i =
1}^{n}K\left(\frac{X_{i\Delta_{n}} - x}{h_{n}}\right) [
M^{2}(X_{i\Delta_{n}}) - M^{2}(x)] }{ \frac{\Delta_{n}}{h_{n}}
\sum_{i = 1}^{n}K\left(\frac{X_{i\Delta_{n}} - x}{h_{n}}\right) }\\
& = & \frac{\frac{1}{h_{n}} \int_{- \infty}^{+ \infty}K\left(\frac{a
- x}{h_{n}}\right)(M^{2}(a) - M^{2}(x)) s(a) da + o_{a.s.}(1)
}{\frac{1}{h_{n}} \int_{- \infty}^{+ \infty}K\left(\frac{a -
x}{h_{n}}\right) s(a) da + o_{a.s.}(1)} + o_{a.s.}(1)\\
& = & h_{n}^{2} K_{2} \left[\frac{1}{2}(M^{2})^{''}(x) +
(M^{2})^{'}(x)\frac{s^{'}(x)}{s(x)} \right] + o_{a.s.}(h_{n}^{2}),
\end{eqnarray*}
with $K_{2} = \int_{- \infty}^{+ \infty} u^{2}K(u)du.$

For the term $V,$ write $\bar{M}^{2}(X_{i\Delta_{n}}) :=
\frac{1}{m_{n}(i\Delta_{n}) \Delta_{n}} \sum_{j =
0}^{m_{n}(i\Delta_{n}) - 1}\left[X_{t(i\Delta_{n})_{j} + \Delta_{n}}
- X_{t(i\Delta_{n})_{j}}\right]^{2},$ we have
\begin{eqnarray*}
V & = & \frac{ \frac{\Delta_{n}}{h_{n}} \sum_{i =
1}^{n}K\left(\frac{X_{i\Delta_{n}} - x}{h_{n}}\right) [
\bar{M}^{2}(X_{i\Delta_{n}}) - M^{2}(X_{i\Delta_{n}})] }{
\frac{\Delta_{n}}{h_{n}}
\sum_{i = 1}^{n}K\left(\frac{X_{i\Delta_{n}} - x}{h_{n}}\right) }\\
& := & \frac{V^{Num}}{V_{Den}}.
\end{eqnarray*}
Due to lemma \ref{lk}, we can obtain that
\begin{equation}
V_{Den} \stackrel{\mathrm{a.s.}}{\longrightarrow}
\bar{L}^{\oplus}_{X}(T, x),
\end{equation}
so we should deal with the term $V^{Num}$ in what follows.
\begin{eqnarray*}
& ~ & V^{Num}\\ & = & \frac{\Delta_{n}}{h_{n}} \sum_{i =
1}^{n}K\left(\frac{X_{i\Delta_{n}} - x}{h_{n}}\right) [
\bar{M}^{2}(X_{i\Delta_{n}}) - M^{2}(X_{i\Delta_{n}})]\\
& = & \frac{\Delta_{n}}{h_{n}} \sum_{i =
1}^{n}K\left(\frac{X_{i\Delta_{n}} - x}{h_{n}}\right) \frac{\sum_{j
= 1}^{n - 1}1_{\{|X_{j\Delta_{n}} - X_{i\Delta_{n}}| \leq
\varepsilon_{n}\}} \frac{1}{\Delta_{n}} \left[(X_{(j+1)\Delta_{n}} -
X_{j\Delta_{n}})^{2} - \Delta_{n}M^{2}(X_{i\Delta_{n}}) \right] }{
\sum_{j = 1}^{n}1_{\{|X_{j\Delta_{n}} - X_{i\Delta_{n}}| \leq
\varepsilon_{n}\}}}\\
& ~ & - \frac{\Delta_{n}}{h_{n}} \sum_{i =
1}^{n}K\left(\frac{X_{i\Delta_{n}} - x}{h_{n}}\right) \frac{
M^{2}(X_{i\Delta_{n}})\frac{\Delta_{n}}{2\varepsilon_{n}}
1_{\{|X_{n\Delta_{n}} - X_{i\Delta_{n}}| \leq \varepsilon_{n}\}} }{
\frac{\Delta_{n}}{2\varepsilon_{n}} \sum_{j =
1}^{n}1_{\{|X_{j\Delta_{n}} - X_{i\Delta_{n}}| \leq
\varepsilon_{n}\}}}\\
& := & V^{Num}_{1} + V^{Num}_{2}.
\end{eqnarray*}
Using the occupation time formula in lemma \ref{lot}, we can
conclude that
\begin{equation}
V^{Num}_{2} =
O_{a.s.}\left(\frac{\Delta_{n}}{\varepsilon_{n}}\right).
\end{equation}
Using It\^{o} formula to the jump-diffusion setting shown in Protter
\cite{pr}, we can write
\begin{eqnarray*}
& ~ & (X_{(j + 1)\Delta_{n}} - X_{j\Delta_{n}})^{2}\\
 & = & 2\int_{j\Delta_{n}}^{(j + 1)\Delta_{n}}(X_{s-} -
 X_{j\Delta_{n}})\mu(X_{s-})ds + 2\int_{j\Delta_{n}}^{(j + 1)\Delta_{n}}(X_{s-} -
 X_{j\Delta_{n}})\sigma(X_{s-})dW_{s}\\
 & ~ & + 2\int_{j\Delta_{n}}^{(j + 1)\Delta_{n}}(X_{s-} -
 X_{j\Delta_{n}}) \int_{Y}c(X_{s-}, y) \bar{\nu}(ds, dy) +
 \int_{j\Delta_{n}}^{(j + 1)\Delta_{n}} M^{2}(X_{s-})ds\\
 & ~ & +
 \int_{j\Delta_{n}}^{(j + 1)\Delta_{n}} \int_{Y}c^{2}(X_{s-}, y) \bar{\nu}(ds,
 dy),
\end{eqnarray*}
which implies that the term $V^{Num}_{1}$ can be divided into five
parts as
\begin{equation}
V^{Num}_{1} = V^{Num}_{11} + V^{Num}_{12} + V^{Num}_{13} +
V^{Num}_{14} + V^{Num}_{15},
\end{equation}
where
\begin{align*}
& V^{Num}_{11} = \frac{\Delta_{n}}{h_{n}} \sum_{i =
1}^{n}K\left(\frac{X_{i\Delta_{n}} - x}{h_{n}}\right) \frac{\sum_{j
= 1}^{n - 1}1_{\{|X_{j\Delta_{n}} - X_{i\Delta_{n}}| \leq
\varepsilon_{n}\}} \frac{1}{\Delta_{n}}
\int_{j\Delta_{n}}^{(j+1)\Delta_{n}}(M^{2}(X_{s-}) -
M^{2}(X_{i\Delta_{n}}))ds }{ \sum_{j = 1}^{n}1_{\{|X_{j\Delta_{n}} -
X_{i\Delta_{n}}| \leq
\varepsilon_{n}\}}},\\
& V^{Num}_{12} = \frac{\Delta_{n}}{h_{n}} \sum_{i =
1}^{n}K\left(\frac{X_{i\Delta_{n}} - x}{h_{n}}\right) \frac{\sum_{j
= 1}^{n - 1}1_{\{|X_{j\Delta_{n}} - X_{i\Delta_{n}}| \leq
\varepsilon_{n}\}} \frac{2}{\Delta_{n}}\int_{j\Delta_{n}}^{(j +
1)\Delta_{n}}(X_{s-} -
 X_{j\Delta_{n}})\mu(X_{s-})ds }{ \sum_{j = 1}^{n}1_{\{|X_{j\Delta_{n}} -
X_{i\Delta_{n}}| \leq
\varepsilon_{n}\}}},\\
& V^{Num}_{13} = \frac{\Delta_{n}}{h_{n}} \sum_{i =
1}^{n}K\left(\frac{X_{i\Delta_{n}} - x}{h_{n}}\right) \frac{\sum_{j
= 1}^{n - 1}1_{\{|X_{j\Delta_{n}} - X_{i\Delta_{n}}| \leq
\varepsilon_{n}\}} \frac{2}{\Delta_{n}}\int_{j\Delta_{n}}^{(j +
1)\Delta_{n}}(X_{s-} -
 X_{j\Delta_{n}})\sigma(X_{s-})dW_{s} }{ \sum_{j = 1}^{n}1_{\{|X_{j\Delta_{n}} -
X_{i\Delta_{n}}| \leq
\varepsilon_{n}\}}},\\
& V^{Num}_{14} = \frac{\Delta_{n}}{h_{n}} \sum_{i =
1}^{n}K\left(\frac{X_{i\Delta_{n}} - x}{h_{n}}\right) \frac{\sum_{j
= 1}^{n - 1}1_{\{|X_{j\Delta_{n}} - X_{i\Delta_{n}}| \leq
\varepsilon_{n}\}} \frac{2}{\Delta_{n}}\int_{j\Delta_{n}}^{(j +
1)\Delta_{n}}(X_{s-} -
 X_{j\Delta_{n}})\int_{Y}c(X_{s-}, y) \bar{\nu}(ds, dy) }{ \sum_{j = 1}^{n}1_{\{|X_{j\Delta_{n}} -
X_{i\Delta_{n}}| \leq
\varepsilon_{n}\}}}, \\
& V^{Num}_{15} = \frac{\Delta_{n}}{h_{n}} \sum_{i =
1}^{n}K\left(\frac{X_{i\Delta_{n}} - x}{h_{n}}\right)
\frac{\frac{1}{2\varepsilon_{n}} \sum_{j = 1}^{n -
1}1_{\{|X_{j\Delta_{n}} - X_{i\Delta_{n}}| \leq \varepsilon_{n}\}}
\int_{j\Delta_{n}}^{(j+1)\Delta_{n}}\int_{Y}c^{2}(X_{s-}, y)
\bar{\nu}(ds, dy) }{\frac{\Delta_{n}}{2\varepsilon_{n}} \sum_{j =
1}^{n}1_{\{|X_{j\Delta_{n}} - X_{i\Delta_{n}}| \leq
\varepsilon_{n}\}}}.
\end{align*}
With the equation (\ref{holder}), we can easily get
\begin{align*}
&V^{Num}_{12} = O_{a.s.}\left(\bar{L}^{\oplus}_{X}(T, x)
\left(\Delta_{n} \log
\left(\frac{1}{\Delta_{n}}\right)\right)^{1/2}\right),\\
 & V^{Num}_{13} = \left(\Delta_{n} \log
\left(\frac{1}{\Delta_{n}}\right)\right)^{1/2}
 O_{P}\left(V^{Num}_{15}\right) = o_{P}\left(V^{Num}_{15}\right),\\
& V^{Num}_{14} = \left(\Delta_{n} \log
\left(\frac{1}{\Delta_{n}}\right)\right)^{1/2}
 O_{P}\left(V^{Num}_{15}\right) = o_{P}\left(V^{Num}_{15}\right).
\end{align*}
For the bias effect term $V^{Num}_{11},$ we have
\begin{eqnarray*}
& ~ & \frac{V^{Num}_{11}}{\frac{\Delta_{n}}{h_{n}} \sum_{i =
1}^{n}K\left(\frac{X_{i\Delta_{n}} - x}{h_{n}}\right)}\\
& = & \frac{\frac{\Delta_{n}}{h_{n}} \sum_{i =
1}^{n}K\left(\frac{X_{i\Delta_{n}} - x}{h_{n}}\right) \frac{\sum_{j
= 1}^{n - 1}1_{\{|X_{j\Delta_{n}} - X_{i\Delta_{n}}| \leq
\varepsilon_{n}\}} \frac{1}{\Delta_{n}}
\int_{j\Delta_{n}}^{(j+1)\Delta_{n}}(M^{2}(X_{j\Delta_{n}}) -
M^{2}(X_{i\Delta_{n}}))ds }{ \sum_{j = 1}^{n}1_{\{|X_{j\Delta_{n}} -
X_{i\Delta_{n}}| \leq \varepsilon_{n}\}}}}{\frac{\Delta_{n}}{h_{n}}
\sum_{i = 1}^{n}K\left(\frac{X_{i\Delta_{n}} - x}{h_{n}}\right)}\\
& ~ & + \frac{\frac{\Delta_{n}}{h_{n}} \sum_{i =
1}^{n}K\left(\frac{X_{i\Delta_{n}} - x}{h_{n}}\right) \frac{\sum_{j
= 1}^{n - 1}1_{\{|X_{j\Delta_{n}} - X_{i\Delta_{n}}| \leq
\varepsilon_{n}\}} \frac{1}{\Delta_{n}}
\int_{j\Delta_{n}}^{(j+1)\Delta_{n}}(M^{2}(X_{s-}) -
M^{2}(X_{j\Delta_{n}}))ds }{ \sum_{j = 1}^{n}1_{\{|X_{j\Delta_{n}} -
X_{i\Delta_{n}}| \leq \varepsilon_{n}\}}}}{\frac{\Delta_{n}}{h_{n}}
\sum_{i = 1}^{n}K\left(\frac{X_{i\Delta_{n}} - x}{h_{n}}\right)}\\
& = & D_{1n} + D_{2n}.
\end{eqnarray*}
By the mean-value theorem, the equation (\ref{holder}) and the
locally boundedness of $(M^{2})^{'}(\cdot),$ it can be shown that
\begin{equation}
D_{2n} = O_{a.s.}\left( \left(\Delta_{n} \log
\left(\frac{1}{\Delta_{n}}\right)\right)^{1/2}\right)
\end{equation}
Moreover,
\begin{eqnarray*}
& ~ & D_{1n}\\ & = & \frac{\frac{1}{h_{n}} \int_{- \infty}^{+
\infty}K\left(\frac{a - x}{h_{n}}\right) \frac{\int_{- \infty}^{+
\infty} 1_{\{|b - a| \leq \varepsilon_{n}\}}(M^{2}(b) -
M^{2}(a))s(b)db}{\int_{- \infty}^{+ \infty} 1_{\{|b - a| \leq
\varepsilon_{n}\}}s(b)db} s(a) da  + o_{a.s.}(1)}{\frac{1}{h_{n}}
\int_{- \infty}^{+ \infty}K\left(\frac{a - x}{h_{n}}\right) s(a) da
+ o_{a.s.}(1) } \\
& \stackrel{\frac{a - x}{h_{n}} = c}{=} & \frac{\frac{1}{h_{n}}
\int_{- \infty}^{+ \infty}K\left(\frac{a - x}{h_{n}}\right)
\frac{\int_{- \infty}^{+ \infty} 1_{\{|b - a| \leq
\varepsilon_{n}\}}(M^{2}(b) - M^{2}(a))s(b)db}{\int_{- \infty}^{+
\infty} 1_{\{|b - a| \leq \varepsilon_{n}\}}s(b)db} s(a)
da}{\frac{1}{h_{n}} \int_{- \infty}^{+ \infty}K\left(\frac{a -
x}{h_{n}}\right) s(a) da} \\
& = & \frac{ \int_{- \infty}^{+ \infty}K\left(c\right) \frac{\int_{-
\infty}^{+ \infty} 1_{\{|\frac{b - x - ch_{n}}{\varepsilon_{n}}|
\leq 1\}}(M^{2}(b) - M^{2}(x + ch_{n}))s(b)db}{\int_{- \infty}^{+
\infty} 1_{\{|\frac{b - x - ch_{n}}{\varepsilon_{n}}| \leq
1\}}s(b)db} s(x + ch_{n}) dc}{ \int_{- \infty}^{+
\infty}K\left(c\right) s(x + ch_{n}) dc} \\
& \stackrel{\frac{b - x}{\varepsilon_{n}} = a}{=} & \frac{ \int_{-
\infty}^{+ \infty}K\left(c\right) \frac{\int_{- \infty}^{+ \infty}
1_{\{|a - \frac{ch_{n}}{\varepsilon_{n}}| \leq 1\}}(M^{2}(x +
a\varepsilon_{n}) - M^{2}(x + ch_{n}))s(x +
a\varepsilon_{n})da}{\int_{- \infty}^{+ \infty} 1_{\{|a -
\frac{ch_{n}}{\varepsilon_{n}}| \leq 1\}}s(x + a\varepsilon_{n})da}
s(x + ch_{n}) dc}{ \int_{- \infty}^{+ \infty}K\left(c\right) s(x +
ch_{n}) dc}\\
& = & \frac{ \int_{- \infty}^{+ \infty}K\left(c\right) \frac{\int_{-
\infty}^{+ \infty} 1_{\{|a - \frac{ch_{n}}{\varepsilon_{n}}| \leq
1\}}(M^{2}(x + a\varepsilon_{n}) - M^{2}(x))s(x +
a\varepsilon_{n})da}{\int_{- \infty}^{+ \infty} 1_{\{|a -
\frac{ch_{n}}{\varepsilon_{n}}| \leq 1\}}s(x + a\varepsilon_{n})da}
s(x + ch_{n}) dc}{ \int_{- \infty}^{+ \infty}K\left(c\right) s(x +
ch_{n}) dc}\\
& ~ & + \frac{ \int_{- \infty}^{+ \infty}K\left(c\right)
\frac{\int_{- \infty}^{+ \infty} 1_{\{|a -
\frac{ch_{n}}{\varepsilon_{n}}| \leq 1\}}(M^{2}(x) - M^{2}(x +
ch_{n}))s(x + a\varepsilon_{n})da}{\int_{- \infty}^{+ \infty}
1_{\{|a - \frac{ch_{n}}{\varepsilon_{n}}| \leq 1\}}s(x +
a\varepsilon_{n})da} s(x + ch_{n}) dc}{ \int_{- \infty}^{+
\infty}K\left(c\right) s(x + ch_{n}) dc}.
\end{eqnarray*}
By use of Taylor expansion, we can formulate
\begin{align*}
& M^{2}(x + a\varepsilon_{n}) -
M^{2}(x) = (M^{2})^{'}(x)a\varepsilon_{n} + \frac{1}{2}(M^{2})^{''}(x)(a\varepsilon_{n})^{2} + o(\varepsilon_{n}^{2}),\\
& M^{2}(x + ch_{n}) - M^{2}(x) = (M^{2})^{'}(x)ch_{n} +
\frac{1}{2}(M^{2})^{''}(x)(ch_{n})^{2} + o(h_{n}^{2}),\\
& s(x + a\varepsilon_{n}) = s(x) + s^{'}(x)a\varepsilon_{n} +
\frac{1}{2}(s)^{''}(x)(a\varepsilon_{n})^{2} + o(\varepsilon_{n}^{2}),\\
& s(x + ch_{n}) = s(x) + s^{'}(x)ch_{n} +
\frac{1}{2}(s)^{''}(x)(ch_{n})^{2} + o(h_{n}^{2}).
\end{align*}

If $h_{n} = o(\varepsilon_{n}),$ we have
\begin{equation}
D_{1n} = \frac{\varepsilon_{n}^{2}}{3}
\left[\frac{1}{2}(M^{2})^{''}(x) +
(M^{2})^{'}(x)\frac{s^{'}(x)}{s(x)}\right] + o(\varepsilon_{n}^{2}).
\end{equation}

If $h_{n} = O(\varepsilon_{n})$ with $\frac{h_{n}}{\varepsilon_{n}}
= \phi,$ then we can get
\begin{eqnarray*}
& ~ & D_{1n}\\ & = & \frac{ \int_{- \infty}^{+
\infty}K\left(c\right) \frac{\int_{- \infty}^{+ \infty} 1_{\{|a -
c\phi| \leq 1\}}((M^{2})^{'}(x)a\varepsilon_{n} +
\frac{1}{2}(M^{2})^{''}(x)(a\varepsilon_{n})^{2})s(x +
a\varepsilon_{n})da}{\int_{- \infty}^{+ \infty} 1_{\{|a - c\phi|
\leq 1\}}s(x + a\varepsilon_{n})da} s(x + ch_{n}) dc}{ \int_{-
\infty}^{+ \infty}K\left(c\right) s(x +
ch_{n}) dc}\\
& ~ & - \frac{ \int_{- \infty}^{+ \infty}K\left(c\right)
\frac{\int_{- \infty}^{+ \infty} 1_{\{|a - c\phi| \leq
1\}}((M^{2})^{'}(x)ch_{n} +
\frac{1}{2}(M^{2})^{''}(x)(ch_{n})^{2})s(x +
a\varepsilon_{n})da}{\int_{- \infty}^{+ \infty} 1_{\{|a - c\phi|
\leq 1\}}s(x + a\varepsilon_{n})da} s(x + ch_{n}) dc}{ \int_{-
\infty}^{+ \infty}K\left(c\right) s(x + ch_{n})
dc}\\
& = & D_{11n} - D_{12n}.
\end{eqnarray*}
Set $g = a- c\phi$, it follows for $D_{11n}$
\begin{eqnarray*}
& ~ & D_{11n}\\
& = & \frac{ \int_{- \infty}^{+ \infty}K\left(c\right) \frac{\int_{-
\infty}^{+ \infty} 1_{\{|g| \leq 1\}}((M^{2})^{'}(x)(g +
c\phi)\varepsilon_{n} + \frac{1}{2}(M^{2})^{''}(x)((g +
c\phi)\varepsilon_{n})^{2})s(x + (g +
c\phi)\varepsilon_{n})dg}{\int_{- \infty}^{+ \infty} 1_{\{|g| \leq
1\}}s(x + (g + c\phi)\varepsilon_{n})dg} s(x + ch_{n}) dc}{ \int_{-
\infty}^{+ \infty}K\left(c\right) s(x + ch_{n}) dc}\\
& = & \frac{ \int_{- \infty}^{+ \infty}K\left(c\right) \frac{\int_{-
\infty}^{+ \infty} 1_{\{|g| \leq 1\}}((M^{2})^{'}(x) g
\varepsilon_{n} + \frac{1}{2}(M^{2})^{''}(x)g^{2}\varepsilon_{n}^{2}
+ (M^{2})^{''}(x) \phi gc\varepsilon_{n}^{2})s(x + (g +
c\phi)\varepsilon_{n})dg}{\int_{- \infty}^{+ \infty} 1_{\{|g| \leq
1\}}s(x + (g + c\phi)\varepsilon_{n})dg} s(x + ch_{n}) dc}{ \int_{-
\infty}^{+ \infty}K\left(c\right) s(x + ch_{n}) dc}\\
& ~ & + h_{n}^{2} K_{2} \left[\frac{1}{2}(M^{2})^{''}(x) +
(M^{2})^{'}(x)\frac{s^{'}(x)}{s(x)}\right] + o(h_{n}^{2})\\
& = & \frac{\frac{1}{2}\int_{- \infty}^{+ \infty} 1_{\{|g| \leq
1\}}g^{2}\varepsilon_{n}^{2}(\frac{1}{2}(M^{2})^{''}(x) s(x) +
(M^{2})^{'}(x) s^{'}(x)) }{s(x)} + o(\varepsilon_{n}^{2})\\
& ~ & + h_{n}^{2} K_{2} \left[\frac{1}{2}(M^{2})^{''}(x) +
(M^{2})^{'}(x)\frac{s^{'}(x)}{s(x)}\right] + o(h_{n}^{2})\\
& = & \frac{1}{3} \varepsilon_{n}^{2}
\left[\frac{1}{2}(M^{2})^{''}(x) +
(M^{2})^{'}(x)\frac{s^{'}(x)}{s(x)}\right] + h_{n}^{2} K_{2}
\left[\frac{1}{2}(M^{2})^{''}(x) +
(M^{2})^{'}(x)\frac{s^{'}(x)}{s(x)}\right] + o(\varepsilon_{n}^{2} +
h_{n}^{2})\\
& = & \varepsilon_{n}^{2}\left(K_{2}\phi^{2} +
\frac{1}{3}\right)\left[\frac{1}{2}(M^{2})^{''}(x) +
(M^{2})^{'}(x)\frac{s^{'}(x)}{s(x)}\right] + o(\varepsilon_{n}^{2}).
\end{eqnarray*}
As for $D_{12n},$ it can be concluded that
\begin{equation}
D_{12n} = h_{n}^{2} K_{2} \left[\frac{1}{2}(M^{2})^{''}(x) +
(M^{2})^{'}(x)\frac{s^{'}(x)}{s(x)} \right] + o_{a.s.}(h_{n}^{2}),
\end{equation}
based on the following similar definite integration
\begin{eqnarray*}
\frac{1}{2} \int_{- \infty}^{+ \infty} 1_{\{|a - c\phi| \leq 1\}}
\frac{ac}{\phi} da = \frac{c}{\phi} \int_{- 1 + c\phi}^{1 + c\phi} a
da = c^{2}.
\end{eqnarray*}
To conclude, when $h_{n} = o(\varepsilon_{n}),$ the bias term for
$\hat{M}^{2}_{n}(x) - M^{2}(x)$ is
\begin{eqnarray}
& ~ & \frac{V^{Num}_{11}}{\frac{\Delta_{n}}{h_{n}} \sum_{i =
1}^{n}K\left(\frac{X_{i\Delta_{n}} - x}{h_{n}}\right)} + B \nonumber\\
& = & D_{11n} - D_{12n} + D_{2n} + B \nonumber\\
& = & \frac{\varepsilon_{n}^{2}}{3} \left[\frac{1}{2}(M^{2})^{''}(x)
+ (M^{2})^{'}(x)\frac{s^{'}(x)}{s(x)} \right] +
o_{a.s.}(\varepsilon_{n}^{2}).
\end{eqnarray}
If $h_{n} = O(\varepsilon_{n})$ with $\frac{h_{n}}{\varepsilon_{n}}
= \phi,$ the total bias term is
\begin{eqnarray}
& ~ & \frac{V^{Num}_{11}}{\frac{\Delta_{n}}{h_{n}} \sum_{i =
1}^{n}K\left(\frac{X_{i\Delta_{n}} - x}{h_{n}}\right)} + B \nonumber\\
& = & D_{11n} - D_{12n} + D_{2n} + B \nonumber\\
& = & \varepsilon_{n}^{2}\left(K_{2}\phi^{2} +
\frac{1}{3}\right)\left[\frac{1}{2}(M^{2})^{''}(x) +
(M^{2})^{'}(x)\frac{s^{'}(x)}{s(x)}\right] + o(\varepsilon_{n}^{2}).
\end{eqnarray}
For the variance effect term $V^{Num}_{15},$ we have
\begin{eqnarray*}
& ~ & \sqrt{\varepsilon_{n}} V^{Num}_{15}\\
& = & \frac{\Delta_{n}}{h_{n}} \sum_{i =
1}^{n}K\left(\frac{X_{i\Delta_{n}} - x}{h_{n}}\right)
\frac{\frac{1}{2 \sqrt{\varepsilon_{n}}} \sum_{j = 1}^{n -
1}1_{\{|X_{j\Delta_{n}} - X_{i\Delta_{n}}| \leq \varepsilon_{n}\}}
\int_{j\Delta_{n}}^{(j+1)\Delta_{n}}\int_{Y}c^{2}(X_{s-}, y)
\bar{\nu}(ds, dy) }{\frac{\Delta_{n}}{2\varepsilon_{n}} \sum_{j =
1}^{n}1_{\{|X_{j\Delta_{n}} - X_{i\Delta_{n}}| \leq
\varepsilon_{n}\}}}.
\end{eqnarray*}
Denote $J_{j} :=
\int_{j\Delta_{n}}^{(j+1)\Delta_{n}}\int_{Y}c^{2}(X_{s-}, y)
\bar{\nu}(ds, dy),$ which is a martingale difference series. By
Gaussian approximation of locally square-integrable martingales
(more technical details seen in Lin and Wang \cite{lw} and Philipp
and Stout \cite{ps}), on an extension of the filtered probability
space we have
\begin{equation}
\label {ga} \max_{1 \leq j \leq n}\left|J_{j} -
B_{\int_{j\Delta_{n}}^{(j+1)\Delta_{n}} \lambda(X_{s-})
\int_{Y}c^{4}(X_{s-}, y)\Pi(dy)ds} \right| = o_{a.s.}(1).
\end{equation}
Based on the equation (\ref{ga}), $\sqrt{\varepsilon_{n}}
V^{Num}_{15}$ has the same asymptotic distribution with
\begin{eqnarray*}
& ~ & E_{n}\\
& := & \frac{\Delta_{n}}{h_{n}} \sum_{i =
1}^{n}K\left(\frac{X_{i\Delta_{n}} - x}{h_{n}}\right)
\frac{\frac{1}{2 \sqrt{\varepsilon_{n}}} \sum_{j = 1}^{n -
1}1_{\{|X_{j\Delta_{n}} - X_{i\Delta_{n}}| \leq \varepsilon_{n}\}}
B_{\int_{j\Delta_{n}}^{(j+1)\Delta_{n}} \lambda(X_{s-})
\int_{Y}c^{4}(X_{s-}, y)\Pi(dy)ds}
}{\frac{\Delta_{n}}{2\varepsilon_{n}} \sum_{j =
1}^{n}1_{\{|X_{j\Delta_{n}} - X_{i\Delta_{n}}| \leq
\varepsilon_{n}\}}},
\end{eqnarray*}
which can be embedded in a time-changed Brownian motion with the
quadratic variation process under assumption in what follows,
\begin{eqnarray*}
& ~ & [E_{n}]\\
& = & \frac{\Delta^{2}_{n}}{h^{2}_{n}} \sum_{i = 1}^{n} \sum_{k =
1}^{n} K\left(\frac{X_{i\Delta_{n}} - x}{h_{n}}\right)
K\left(\frac{X_{k\Delta_{n}} - x}{h_{n}}\right) \times\\
& \times & \frac{\frac{1}{4 \varepsilon_{n}} \sum_{j = 1}^{n -
1}1_{\{|X_{j\Delta_{n}} - X_{i\Delta_{n}}| \leq \varepsilon_{n}\}}
1_{\{|X_{j\Delta_{n}} - X_{k\Delta_{n}}| \leq \varepsilon_{n}\}}
\int_{j\Delta_{n}}^{(j+1)\Delta_{n}} \lambda(X_{s-})
\int_{Y}c^{4}(X_{s-}, y)\Pi(dy)ds
}{\left(\frac{\Delta_{n}}{2\varepsilon_{n}} \sum_{j =
1}^{n}1_{\{|X_{j\Delta_{n}} - X_{i\Delta_{n}}| \leq
\varepsilon_{n}\}}\right) \left(\frac{\Delta_{n}}{2\varepsilon_{n}}
\sum_{j = 1}^{n}1_{\{|X_{j\Delta_{n}} - X_{k\Delta_{n}}| \leq
\varepsilon_{n}\}}\right)}\\
& = & \frac{1}{h_{n}^{2}} \int_{0}^{T} ds \int_{0}^{T} du
K\left(\frac{X_{s-} - x}{h_{n}}\right) K\left(\frac{X_{u-} -
x}{h_{n}}\right) \times\\
& \times & \frac{\frac{1}{4 \varepsilon_{n}} \int_{0}^{T} db
1_{\{|X_{b-} - X_{s-}| \leq \varepsilon_{n}\}} 1_{\{|X_{b-} -
X_{u-}| \leq \varepsilon_{n}\}} M^{4}(X_{b-} +
o_{a.s.}(1))}{\left(\frac{1}{2\varepsilon_{n}} \int_{0}^{T}
1_{\{|X_{b-} - X_{s-}| \leq \varepsilon_{n}\}} db \right)
\left(\frac{1}{2\varepsilon_{n}} \int_{0}^{T} 1_{\{|X_{b-}
- X_{u-}| \leq \varepsilon_{n}\}} db \right)}\\
& = & \frac{1}{h_{n}^{2}} \int_{- \infty}^{+ \infty} ds \int_{-
\infty}^{+ \infty} du K\left(\frac{s - x}{h_{n}}\right)
K\left(\frac{u -
x}{h_{n}}\right) \times\\
& \times & \frac{\frac{1}{4 \varepsilon_{n}} \int_{- \infty}^{+
\infty} db 1_{\{|b - s| \leq \varepsilon_{n}\}} 1_{\{|b - u| \leq
\varepsilon_{n}\}} M^{4}(b) \bar{L}_{X}(T, b) \bar{L}_{X}(T, s)
\bar{L}_{X}(T, u)}{\left(\frac{1}{2\varepsilon_{n}} \int_{-
\infty}^{+ \infty} 1_{\{|b - s| \leq \varepsilon_{n}\}}
\bar{L}_{X}(T, b) db \right) \left(\frac{1}{2\varepsilon_{n}}
\int_{- \infty}^{+ \infty} 1_{\{|b - u| \leq \varepsilon_{n}\}}
\bar{L}_{X}(T,
b) db \right)} + o_{a.s.}(1)\\
& \stackrel{\frac{s - x}{h_{n}} = a~\frac{u - x}{h_{n}} = e}{=} &
o_{a.s.}(1) + \frac{1}{4 \varepsilon_{n}} \int_{- \infty}^{+ \infty}
da \int_{-
\infty}^{+ \infty} de K\left(a\right) K\left(e\right) \times\\
& \times & \frac{ \int_{- \infty}^{+ \infty} db 1_{\{|b - x -
ah_{n}| \leq \varepsilon_{n}\}} 1_{\{|b - x - eh_{n}| \leq
\varepsilon_{n}\}} M^{4}(b) \bar{L}_{X}(T, b) \bar{L}_{X}(T, x +
ah_{n}) \bar{L}_{X}(T, x + eh_{n})}{\left(\frac{1}{2\varepsilon_{n}}
\int_{- \infty}^{+ \infty} 1_{\{|b - x - ah_{n}| \leq
\varepsilon_{n}\}} \bar{L}_{X}(T, b) db \right)
\left(\frac{1}{2\varepsilon_{n}} \int_{- \infty}^{+ \infty} 1_{\{|b
- x - eh_{n}| \leq \varepsilon_{n}\}} \bar{L}_{X}(T, b) db
\right)}\\
& \stackrel{\frac{b - x}{\varepsilon_{n}} = z }{=} & o_{a.s.}(1) +
\frac{1}{4 } \int_{- \infty}^{+ \infty} da \int_{-
\infty}^{+ \infty} de K\left(a\right) K\left(e\right) \times\\
& \times & \frac{ \int_{- \infty}^{+ \infty} dz 1_{\{|z -
a\frac{h_{n}}{\varepsilon_{n}}| \leq 1\}} 1_{\{|z -
e\frac{h_{n}}{\varepsilon_{n}}| \leq 1\}} M^{4}(x +
z\varepsilon_{n}) \bar{L}_{X}(T, x + z\varepsilon_{n})
\bar{L}_{X}(T, x + ah_{n}) \bar{L}_{X}(T, x +
eh_{n})}{\left(\frac{1}{2} \int_{- \infty}^{+ \infty} 1_{\{|z -
a\frac{h_{n}}{\varepsilon_{n}}| \leq 1\}} \bar{L}_{X}(T, x +
z\varepsilon_{n}) dz \right) \left(\frac{1}{2} \int_{- \infty}^{+
\infty} 1_{\{|z - e\frac{h_{n}}{\varepsilon_{n}}| \leq 1\}}
\bar{L}_{X}(T, x + \varepsilon_{n}) dz \right)}.
\end{eqnarray*}
To conclude, when $h_{n} = o(\varepsilon_{n}),$ we have
\begin{equation}
[E_{n}] \stackrel{a.s.}{\longrightarrow} \frac{1}{2} M^{4}(x)
\bar{L}_{X}(T, x).
\end{equation}
So with lemma \ref{lk} we have
\begin{equation}
\sqrt{\varepsilon_{n} \hat{\bar{L}}_{X}(T, x)}
\frac{V^{Num}_{15}}{\frac{\Delta_{n}}{h_{n}} \sum_{i =
1}^{n}K\left(\frac{X_{i\Delta_{n}} - x}{h_{n}}\right)} \Rightarrow
N\left(0, \frac{1}{2} M^{4}(x)\right),
\end{equation}
which implies with the equation (4.18)
\begin{equation}
\label{r1} \sqrt{\varepsilon_{n} \hat{\bar{L}}_{X}(T, x)}
(\hat{M}^{2}_{n}(x) - M^{2}(x) - \Gamma_{M^{2}}) \Rightarrow
N\left(0, \frac{1}{2} M^{4}(x)\right),
\end{equation}
where $\Gamma_{M^{2}} = \frac{\varepsilon_{n}^{2}}{3}
\left[\frac{1}{2}(M^{2})^{''}(x) +
(M^{2})^{'}(x)\frac{s^{'}(x)}{s(x)} \right].$

If $h_{n} = O(\varepsilon_{n})$ with $\frac{h_{n}}{\varepsilon_{n}}
= \phi,$ we can obtain
\begin{equation}
[E_{n}] \stackrel{a.s.}{\longrightarrow} \frac{1}{2} \theta_{\phi}
M^{4}(x) \bar{L}_{X}(T, x),
\end{equation}
where
\begin{eqnarray*}
\theta_{\phi} & = & \int_{- \infty}^{+ \infty} da \int_{- \infty}^{+
\infty} de K\left(a\right) K\left(e\right) \frac{\frac{1}{2} \int_{-
\infty}^{+ \infty} dz 1_{\{|z - a\phi| \leq 1\}} 1_{\{|z - e\phi|
\leq 1\}} }{\left(\frac{1}{2} \int_{- \infty}^{+ \infty} 1_{\{|z -
a\phi| \leq 1\}} dz \right) \left(\frac{1}{2} \int_{- \infty}^{+
\infty} 1_{\{|z - e\phi| \leq 1\}} dz \right)}\\
& = & \int_{- \infty}^{+ \infty} da \int_{- \infty}^{+ \infty} de
K\left(a\right) K\left(e\right) \frac{1}{2} \int_{- \infty}^{+
\infty} dz 1_{\{|z - a\phi| \leq 1\}} 1_{\{|z - e\phi| \leq 1\}}\\
& = & \frac{1}{2}\int_{- \infty}^{+ \infty} dz \int_{(z -
1)/\phi}^{(z + 1)/\phi} da \int_{(z - 1)/\phi}^{(z + 1)/\phi} de
K\left(a\right) K\left(e\right)
\end{eqnarray*}
So with lemma \ref{lk} we have
\begin{equation}
\sqrt{\varepsilon_{n} \hat{\bar{L}}_{X}(T, x)}
\frac{V^{Num}_{15}}{\frac{\Delta_{n}}{h_{n}} \sum_{i =
1}^{n}K\left(\frac{X_{i\Delta_{n}} - x}{h_{n}}\right)} \Rightarrow
N\left(0, \frac{1}{2} \theta_{\phi} M^{4}(x)\right),
\end{equation}
which implies with the equation (4.19)
\begin{equation}
\label{r2} \sqrt{\varepsilon_{n} \hat{\bar{L}}_{X}(T, x)}
(\hat{M}^{2}_{n}(x) - M^{2}(x) - \Gamma_{M^{2}}^{\phi}) \Rightarrow
N\left(0, \frac{1}{2} \theta_{\phi} M^{4}(x)\right),
\end{equation}
where $\Gamma_{M^{2}}^{\phi} =
\varepsilon_{n}^{2}\left(K_{2}\phi^{2} +
\frac{1}{3}\right)\left[\frac{1}{2}(M^{2})^{''}(x) +
(M^{2})^{'}(x)\frac{s^{'}(x)}{s(x)}\right].$

We have proved the main results in Theorem \ref{mr} based on
(\ref{r1}) and (\ref{r1}).\end{proof}

{\bf Acknowledgments} This research work is supported by the General
Research Fund of Shanghai Normal University (No. SK201720) and
Funding Programs for Youth Teachers of Shanghai Colleges and
Universities (No. A-9103-17-041301).

\bibliographystyle{amsplain}

\end{document}